\documentclass{article}
\usepackage{amsthm}
\usepackage{amsmath}
\usepackage{algorithm}
\usepackage{algpseudocode}
\usepackage{multirow}
\usepackage{amssymb}
\usepackage{indentfirst}
\usepackage{graphicx}
\usepackage{xcolor}
\usepackage{tabularray}
\usepackage[left=1.20in,right=1.20in,top=0.8in,bottom=0.8in]{geometry}
\title{A novel algorithm for optimizing bundle adjustment in image sequence alignment }
\date{}
\author{Hailin Xu, Hongxia Wang, Renmin Han$^*$, Huanshui Zhang$^*$}
\begin{document}
	\maketitle

	\section*{Abstract}
The Bundle Adjustment (BA) model is commonly optimized using a nonlinear least squares method, with the Levenberg-Marquardt (L-M) algorithm being a typical choice. However, despite the L-M algorithm's effectiveness, its sensitivity to initial conditions often results in slower convergence when applied to poorly conditioned datasets, motivating the exploration of alternative optimization strategies. This paper introduces a novel algorithm for optimizing the BA model in the context of image sequence alignment for cryo-electron tomography, utilizing optimal control theory to directly optimize general nonlinear functions. The proposed Optimal Control Algorithm (OCA) exhibits superior convergence rates and effectively mitigates the oscillatory behavior frequently observed in L-M algorithm. Extensive experiments on both synthetic and real-world datasets were conducted to evaluate the algorithm's performance. The results demonstrate that the OCA achieves faster convergence compared to the L-M algorithm. Moreover, the incorporation of a bisection-based update procedure significantly enhances the OCA's performance, particularly in poorly initialized datasets. These findings indicate that the OCA can substantially improve the efficiency of 3D reconstructions in cryo-electron tomography.  \\
	\textbf{Key words:} Bundle adjustment, Electron tomography, Image sequence alignment, Optimal control algorithm.

	\section{Introduction}
	Electron tomography (ET) techniques is an important branch of cryo-electron microscopy. It utilizes the 3D reconstruction method in simultaneous localization and mapping (SLAM) to reconstruct the 3D structure of organisms by photographing biological samples at different oblique angles. 
	\footnotetext{This work was supported by the Original Exploratory Program Project
		of National Natural Science Foundation of China (62450004), the Joint
		Funds of the National Natural Science Foundation of China (U23A20325),
		and the Natural Science Foundation of Shandong
		Province (ZR2021ZD14, ZR2024MF045). (Corresponding author: Renmin Han, Huanshui Zhang.)
		\par H. Zhang is with the College of Electrical Engineering and Automation,
		Shandong University of Science and Technology, Qingdao 266590, China,
		and also with the School of Control Science and Engineering, Shandong
		University, Jinan, Shandong 250061, China (e-mail: hszhang@sdu.edu.cn).
		\par R. Han is with the Research Center for Mathematics and interdisciplinary Sciences, Shandong University, Qingdao 266237, China. (e-mail: hanrenmin@sdu.edu.cn).
		 \par H. Wang is with the College of Electrical Engineering and Automation, Shandong University
		of Science and Technology, Qingdao 266590, China (e-mail: whx1123@126.com;).
		\par  H. Xu is with the College of Mathematics and Systems Science, Shandong University
		of Science and Technology, Qingdao 266590, China (e-mail: 201901070720@sdu.edu.cn;)}This method enables the acquisition of relatively accurate biological structural information. The most critical step is the alignment of image sequence, which has a direct impact on the final 3D structure restoration accuracy. Image sequence alignment can be divided into two main methods: marker-free alignment \cite{GUCKENBERGER1982167} and marker-based alignment \cite{Lawrence1992,autoHAN2015403}. In the marker-based alignment , the bundle adjustment is the core, and the nonlinear least square method is usually used to minimize the reprojection error and optimize the target parameters.

	Up to now, there exists a variety of alternative algorithms for solving nonlinear least squares problems, including the Newton method, Gauss-Newton method, and L-M algorithm, among others. The L-M algorithm is commonly employed in many practical applications. Recently, an optimization algorithm based on optimal control has been proposed \cite{zhangzgkx}. This novel approach transforms the optimization problem into an optimal control problem by designing a cost function closely related to the objective function. An optimal controller adjusts a first-order difference equation to minimize the objective function.
	This method offers several advantages, such as relatively flexible initial value selection, rapid convergence rates, and the avoidance of oscillatory behavior commonly associated with gradient methods. Notably, when the initial value is appropriately chosen, the OCA can yield different local minima by adjusting the input weight matrix.

	At present, the OCA remains primarily at the theoretical level and has not yet seen widespread application. This paper explores the practical application of the OCA by optimizing the BA model in cryo-ET and comparing the results with those obtained using the L-M algorithm. The remainder of this paper is organized as follows. Section 2 provides a review of recent research on BA and offers a brief introduction to the steps of the L-M algorithm. Section 3 discusses the projection model employed in BA within the context of cryo-ET, along with the underlying principles of the OCA. In Section 4, we present a comparative experiment between the L-M algorithm and OCA, demonstrating that the latter exhibits faster convergence. Finally, Section 5 outlines future directions for the development of OCA and explores its potential applications.
    \section{Related work and L-M algorithm}
    	Bundle adjustment refers to the extraction of optimal 3D models and camera parameters (both intrinsic and extrinsic) from visual reconstruction data. It involves adjusting the camera orientation and the spatial positions of feature points such that the optimal alignment of light rays, reflected from each feature point, is achieved. These rays are eventually projected onto the camera's optical center\cite{BAmodel}.

    	The BA was first introduced by Brown and Duane\cite{brown1976bundle}, with its cost function formulated as the sum of squares of multiple nonlinear functions. Given that the minimization process is carried out using nonlinear least squares, the L-M algorithm is the most commonly employed optimization technique. Recent studies have primarily concentrated on expanding the application scenarios of the BA model. For example, Liu et al.\cite{liu2021balm} propose a new cost function that extends the applicability of BA to lidar mapping. However, relatively few works have addressed the fundamental aspects of the problem by exploring alternative algorithmic structures to replace the L-M algorithm. In this regard, Li et al.\cite{articleBFGS} introduce a Broyden–Fletcher–Goldfarb–Shanno algorithm to approximate the cost function. Similarly, Zhou et al. \cite{zhou2020efficient} propose a point-to-plane cost function aimed at jointly optimizing depth camera poses and plane parameters for 3D reconstruction.

    	The aforementioned approaches, including the L-M algorithm, can be regarded as variations of the Gauss-Newton method. Consequently, the linearization of the cost function is an inherent aspect of these methods, which leads to increased errors and a higher number of iterations. In contrast, the OCA, introduced in the subsequent section, does not rely on such linearization techniques.
    	\subsection{Levenberg-Marquardt algorithm}
The L-M algorithm is a widely used algorithm for solving non-linear least squares problems, which combines the concepts of gradient descent and the Gauss-Newton algorithm to provide a robust approach for parameter estimation.

Let $\varphi(x)=(\varphi_1(x),\varphi_2(x),...,\varphi_m(x))^{\top}$, the objective is to minimize the sum of squared residuals
\begin{equation}
	\underset{x \in \mathbb{R}^n}{\text{min}} \quad  \Phi(x)= \frac{1}{2}||\varphi(x) ||^2 =\frac{1}{2} \sum_{i=1}^{m} \varphi_i^2(x) .
\end{equation}

The update rule for the parameters $x$ is given by
\begin{equation}
	\begin{aligned}
		\left(J^{\top}(x_k)J(x_k) + \mu_k I\right)d_k & =-J^{\top}(x_k)\varphi(x_k),\\
		x_{k+1} & =x_k+d_k, 
	\end{aligned}
\end{equation}
where $J(x_k)$ is the Jacobian matrix of partial derivatives, $I$ is the identity matrix and $\mu_k$ is referred to as damping parameter, which controls the behavior of the algorithm.

A larger $\mu_k$ makes the update resemble gradient descent while a smaller $\mu_k$ approximates the Gauss-Newton method. This adaptability allows the L-M algorithm to efficiently handle a wide range of problem conditions.

In this approach, we employ the trust-region method to adjust $\mu_k$. First, we define a quadratic function at the current iteration point to represent the predicted value
\begin{equation}
	P(d)=\Phi(x_k)+(J^{\top}(x_k) \varphi(x_k))^\top d + \frac{1}{2}d^\top(J^{\top}(x_k)J(x_k))d. 
\end{equation}

The step $d_k$ is then computed based on the current value of $\mu_k$, and the ratio of the actual decrease to the predicted decrease is evaluated as
\begin{equation}
	\xi_k=\frac{ \Phi(x_{k+1}) - \Phi(x_{k}) }{ P(d_k) - \Phi(x_{k})}.\label{etak}
\end{equation}

Equation (\ref{etak}) reflects the degree to which the linearized model aligns with the nonlinear objective function.
\begin{itemize}
	\item When $\xi_k$ approaches 1, the quadratic model $P(d)$ provides a good approximation of the objective function at $x_k$, suggesting that  $\mu_k$ should be reduced.
	\item When $\xi_k$ approaches 0, the quadratic model $P(d)$ poorly approximates the objective function at $x_k$, indicating that $\mu_k$ should be  increased.
	\item  When $\xi_k$ is neither close to 0 nor 1, $\mu_k$ is considered to be appropriately chosen, and no adjustment is necessary.
\end{itemize}

Typically, the threshold values for $\xi_k$ are set to 0.25 and 0.75. Based on this, the update rule for $\mu_k$ in the  L-M  method is as follows
\begin{equation}
	\mu_{k+1} : = \left\{
	\begin{aligned}
		10\mu_k, \quad & \xi_k < 0.25 ,\\
		\mu_k,\quad & 0.25 \leq \xi_k \leq 0.75 ,\\
		0.1\mu_k\quad & \xi_k > 0.75 .
	\end{aligned}\right.
\end{equation}

In summary, the complete flow of the L-M algorithm is as follows
\begin{algorithm}
	\caption{Levenberg-Marquardt  Algorithm}
	\begin{algorithmic}[1]
		\State \textbf{Initialize:} Set initial guess $x_0$, initial damping parameter $\mu_0$, and tolerance $\varepsilon$.
		\State \textbf{Input:} Residual function $\varphi(x)$, Jacobian $J(x)$, damping parameter $\mu$.
		\For{$k = 0, 1, 2, \dots$}
		\State Compute the Jacobian $J(x_k)$ and the residual $\varphi(x_k)$.
		\State Solve for $d_k$ using: 
		\[
		\left(J^{\top}(x_k) J(x_k) + \mu_k I\right) d_k = -J^{\top}(x_k) \varphi(x_k)
		\]
		\State Compute the iteration: $x_{k+1} = x_k + d_k$.
		\If{$||x_{k+1}-x_k|| < \varepsilon$}
		\State Terminate the algorithm.
		\EndIf
		\State Evaluate the ratio of actual to update $\mu$:
		\[
		\xi_k = \frac{\varphi(x_k) - \varphi(x_{k+1})}{P(d_k) - \varphi(x_k)}
		\]
		\If{$\xi_k < 0.25$} 
		\State Increase the damping parameter: $\mu_{k+1} = 10\mu_k$.
		\ElsIf{$0.25 \leq \xi_k \leq 0.75$}
		\State Hold the damping parameter: $\mu_{k+1} = \mu_k$.
		\Else
		\State Decrease the damping parameter: $\mu_{k+1} = 0.1\mu_k$.
		\EndIf
		\EndFor
		\State \textbf{Output:} Optimal solution $x^*$.
	\end{algorithmic}
\end{algorithm}

	\section{Method}
This section is divided into two parts. The first part introduces the projection model of cryo-ET, a crucial component in BA. The second part explains the principles of the OCA and extends the algorithm by incorporating a parameter update mechanism.
	\subsection{Projection model}
An essential component of BA is the camera projection model, which maps the 3D coordinates of objects in the real world onto 2D coordinates in an image. In this paper, we focus on the application of a projection model specifically for cryo-ET image sequence\cite{campro}
\begin{equation}
	\begin{pmatrix}
		u \\
		v
	\end{pmatrix}
	= \textbf{R}_{\gamma}^{-1}\left(\frac{1}{s}\textbf{P}\textbf{R}_{\beta}\textbf{R}_{\alpha}
	\begin{pmatrix}
		X\\
		Y\\
		Z
	\end{pmatrix}
	-\textbf{t}\right), \label{op}
\end{equation}
where $(X,Y,Z)^{T}$ is the 3D coordinates of the fiducial markers, $(u,v)^{\top}$ is the 2D projection of the fiducial markers in the image, $s$ is the image scale change, $\textbf{t}=(t_0,t_1)^{\top}$ is the translation during projection and $\textbf{P}$ represents the orthogonal projection matrix. $\textbf{R}$ denotes the rotation matrix of the projection process, the details of $\textbf{R}$ and $\textbf{P}$ as follows
\begin{equation*}
	\textbf{P}=
	\begin{pmatrix}
		1 & 0 & 0 \\
		0 & 1 & 0
	\end{pmatrix},
\end{equation*}
\begin{equation*}
	\textbf{R}_{\alpha}=
	\begin{pmatrix}
		1 & 0 & 0 \\
		0 & cos\alpha & sin\alpha \\
		0 & -sin\alpha & cos\alpha
	\end{pmatrix},
\end{equation*}
\begin{equation*}
	\textbf{R}_{\beta}=
	\begin{pmatrix}
		cos\beta & 0 & -sin\beta \\
		0 & 1 & 0 \\
		sin\beta & 0 & cos\beta
	\end{pmatrix},
\end{equation*}
\begin{equation*}
	\textbf{R}_{\gamma}=
	\begin{pmatrix}
		cos\gamma & sin\gamma \\
		-sin\gamma & cos\gamma
	\end{pmatrix}.
\end{equation*}
$\textbf{x}_c=(s,\alpha,\beta,\gamma,t_0,t_1)^{\top}$ is called as camera parameters and $\textbf{x}_p=(X,Y,Z)^{\top}$ is called as point parameters.

 When sufficient projection data for a single 3D point is available, triangulation can be used to approximate the coordinates of that point. However, due to various uncertainties, the initial estimates of the 3D coordinates may not be entirely accurate. When this estimated 3D point is reprojected onto the 2D image, a discrepancy, or error, arises in comparison to the observed point.

Let the observed point be denoted as $\textbf{z} \stackrel{\Delta}{=} (u,v)^{\top}$ and the point obtained by reprojection can be expressed as follows
\begin{equation}
	\begin{pmatrix}
		\hat{u} \\ \hat{v}
	\end{pmatrix}
=
\textbf{R}_{\gamma_0}^{-1}\left(\frac{1}{s_0}\textbf{P}\textbf{R}_{\beta_0}\textbf{R}_{\alpha_0}
\begin{pmatrix}
	X_0\\
	Y_0\\
	Z_0
\end{pmatrix} 
-\textbf{t}_{0}\right) = \textbf{h}(\textbf{x}_{c_0},\textbf{x}_{p_0}), \label{rep}
\end{equation}
where $(\cdot)_0$ represents the initial estimation of parameters. Thus the reprojection error can be written as 
\begin{equation}
	\textbf{e}=\textbf{z}-\textbf{h}. \label{reperror}
\end{equation}

To account for all observations, we can introduce a subscript to represent the error for each observation. Let $\textbf{z}_{ij}$ denotes the data generated by the $i$-th marker in the $j$-th image. The overall cost function can be expressed as
\begin{equation}
	\frac{1}{2} \sum_{i=1}^{n} \sum_{j=1}^{m} ||\textbf{e}_{ij}||^2 \delta_{ij}
	=
	\sum_{i=1}^{n} \sum_{j=1}^{m}
	\frac{1}{2}\left|\left|\textbf{z}_{ij} - \textbf{h}(\textbf{x}_{c_j},\textbf{x}_{p_i})\right|\right|^2\delta_{ij}.
	\label{BAcost}
\end{equation}

During the sample collection process, some marker points may become obscured at certain angles and thus will not be visible. To represent the visibility of these points, we introduce a visibility indicator $\delta_{ij}=\{0,1\}$. If the projection of the $j$-th marker is visible in the $i$-th image, $\delta_{ij}=1$; otherwise, $\delta_{ij}=0$.

Solving this least-squares problem involves jointly optimizing both the camera parameters $\textbf{x}_c$ and the point parameters $\textbf{x}_p$, a process known as Bundle Adjustment.

The effectiveness of the adjustments is evaluated using the following average residual formula
\begin{equation}
	L_1=\frac{1}{2mn}\sum_{i=1}^{n} \sum_{j=1}^{m}\left|\left|\textbf{z}_{ij} - \textbf{h}(\textbf{x}^*_{c_j},\textbf{x}^*_{p_i})\right|\right|_1,  \label{l1res}
\end{equation}
where $|| \cdot ||_1$ denotes the 1-norm, $x_c^*$ and $x_p^*$ are the adjusted parameters. A smaller $L_1$ value indicates a lower reprojection error.

\subsection{Optimal control algorithm}
Suppose $f(x) : \mathbb{R}^n \rightarrow \mathbb{R}$ is a nonlinear function with a continuous second derivative. Consider the following optimization problem: 
\begin{equation}
	\underset{x \in \mathbb{R}^n}{\text{min}} \quad f(x). \label{OP}
\end{equation}

Generally, it can be solved by the following iteration:
\begin{equation}
	x_{k+1}	=x_k+d_k , \label{ite}
\end{equation}
where $d_k$ is the step size.

If we interpret $d_k$ as a control input in a control system, the iteration (\ref{ite}) exhibits a structure analogous to that of the following discrete-time linear time-invariant system:

\begin{equation}
	x_{k+1}	=x_k+u_k.
\end{equation}

By designing a performance index related to the cost function, the optimization problem (\ref{OP}) can be reformulated as the following optimal control problem:
\begin{equation}
	\begin{aligned}
		& \underset{u}{\text{min}} \quad \sum_{k=1}^{N}(f(x_k)+\frac{1}{2}d_k^\top Rd_k)+f(x_{N+1}), \\
		& s.t. \quad x_{k+1}	=x_k+d_k,
	\end{aligned} \label{OCP}
\end{equation}
where $f(x_{N+1})$ represents the terminal cost function, $N$ is the terminal time, and $x_0$ is given. The positive definite matrix $R$ denotes the control weight.

Based on the Pontryagin Maximum Principle, the optimal trajectory satisfies a set of forward and backward equations. Solving these equations yields the optimal controller for the problem (\ref{OCP}):
\begin{equation}
	d_k=-R^{-1}\sum_{i=k+1}^{N+1}\nabla f(x_i).
\end{equation}

Thus, we can derive the following iterative algorithm:
\begin{equation}
	x_{k+1}=x_k-R^{-1}\sum_{i=k+1}^{N+1}\nabla f(x_i). \label{impite}
\end{equation}

Although the above iteration is implicit, it can be obtained as an approximate explicit solution through Taylor expansion. First, take the first-order Taylor expansion of $\nabla f(x_i)$ at $x_{i-1}$:
\begin{gather}
	x_{k+1}=x_k-d_k(x_k), \quad k=0,...,N, \\
		d_N(x_k)=[R+H(x_k)]^{-1}f(x_N),  \\
	d_l(x_k)=[R+H(x_k)]^{-1}[\nabla f(x_k)+\sum_{i=k+1}^{N}(\nabla f(x_i)-H(x_i)d_{i}(x_k))],  \quad l=N-1, \dots, k \label{taysum}, 
\end{gather}
where $H(x)$ is the Hessian matrix of $f(x)$. By substituting $x_k$ for all $x_i (i > k)$ on the right-hand side of equation (\ref{taysum}) and transforming the backward iteration scheme into a forward formulation, we derive the following optimal control algorithm to solve the BA model:

\begin{gather}
	x_{k+1}=x_k-d_k(x_k),  \quad k=0,...,N, \notag \\
	d_0(x_k)=[R+H(x_k)]^{-1}\nabla f(x_k),  \\ 
	d_l(x_k)=[R+H(x_k)]^{-1}[\nabla f(x_k)+Rd_{l-1}(x_k)], \quad l=1, \dots ,k \notag 
\end{gather}

The optimal control algorithm has been demonstrated to exhibit a superlinear convergence and operates through the following stages. The details can be referred to \cite{zhangzgkx}. 
\begin{algorithm}
	\caption{Optimal Control Algorithm}
	\begin{algorithmic}[1]
		\State \textbf{Initialize:} Set initial guess $x_0$, control weight matrix $R$, and tolerance $\varepsilon$.
		\State \textbf{Input:} Objective function $f(x)$, gradient
		 $\nabla f(x)$, Hessian  $H(x)$.
		\For{$k = 0, 1, 2, \dots$}
		\State Compute the gradient $J(x_k)$ and the Hessian  $H(x_k)$.
		\For{$j = 0, 1, 2, \dots, k$}
		 \If{$j = 0$}
		\State Compute the initial search direction:
		\[
		g_0(x_k) = \left[ R + H(x_k) \right]^{-1} \nabla f(x_k)
		\]
		\Else
		\State Compute the search direction for $j \geq 1$:
		\[
		g_j(x_k) = \left[ R + H(x_k) \right]^{-1} \left[\nabla f(x_k) + R g_{j-1}(x_k)\right]
		\]
		\EndIf
		\EndFor
		\State Update the solution:
		\[
		x_{k+1} = x_k - g_k(x_k)
		\]
		\If{$||x_{k+1}-x_k|| < \varepsilon$}
		\State Terminate the algorithm.
		\EndIf
		\EndFor
		\State \textbf{Output:} Optimal solution $x^*$.
	\end{algorithmic}
\end{algorithm}

 The selection of an appropriate weight matrix $ R $ is critical for the effective operation of the OCA. For convenience, $ R $ is generally defined as $ R = \lambda \cdot I $, where $ \lambda > 0 $ is a constant. Under the condition that $ (R + H(x_0)) > 0 $, smaller values of $ \lambda $ lead to faster convergence of the OCA. In the simulation experiments presented in Section 4.2, a set of poorly initialized data was generated, with initial points positioned far from the extremum points. In such cases, a larger $ \lambda $ is typically required to ensure algorithmic convergence; however, this also results in a slower convergence rate for the OCA. To resolve this issue, we propose an update procedure based on the bisection method, which allows the weight matrix $ R$ to progressively decrease as the number of iterations increases. This approach ensures algorithmic convergence while simultaneously enhancing convergence speed. The detailed update strategy is presented as follows.
\begin{algorithm}
	\caption{Optimal Control Algorithm with Adaptive Update of $R$}
	\begin{algorithmic}[1]
		\State \textbf{Initialize:} Set initial guess $x_0$, update factor $\lambda_0, \lambda_1$, and tolerance $\varepsilon$.
		\State \textbf{Input:} Objective function $f(x)$, gradient
		$\nabla f(x)$, Hessian  $H(x)$.
		\For{$k = 0, 1, 2, \dots$}
		\If{k$\leq$1}
		\State Set $R_k=\lambda_k \cdot I$
		\State Compute the gradient $J(x_k)$ and the Hessian  $H(x_k)$.
		\For{$j = 0, 1, 2, \dots, k$}
		\If{$j = 0$}
		\State Compute the initial search direction:
		\[
		g_0(x_k) = \left[ R_k + H(x_k) \right]^{-1} \nabla f(x_k)
		\]
		\Else
		\State Compute the search direction for $j \geq 1$:
		\[
		g_j(x_k) = \left[ R_k + H(x_k) \right]^{-1} \left[\nabla f(x_k) + R g_{j-1}(x_k)\right]
		\]
		\EndIf
		\EndFor
		\Else
        \State Take $R_k=R_{k-1}$ and do step 7-13.
        \State Set $x_{\text{temp}}=x_k - g_k(x_k)$ and compute the objective function 
        \[cost_1=f(x_{\text{temp}}).\]
        \State Initialize interval bounds for line search: $a = 0$, $b = \lambda_{k-1}$.
        \While{$(b - a) > 0.1$} \Comment{Binary search to adjust $R_k$}
        \State Set $c = \frac{b + a}{2}$ and update $R_k = c \cdot I$
        \State Recompute the search direction using steps 7-13
        \State Set $x_{\text{temp}}=x_k - g_k(x_k)$ and compute the new objective function 
           \[cost_2=f(x_{\text{temp}}).\]
        \If{$cost_1 > cost_2$}
        \State $b = c$
        \State Update $cost_1 = cost_2$
        \ElsIf{$cost_1 < cost_2$}
        \State $a = c$
        \State Update $cost_1 = cost_2$
        \Else
		\State \textbf{break} \Comment{Terminate the binary search if no improvement}
        \EndIf
        \EndWhile
		\State Finalize $R_k = c \cdot I$ \Comment{Set the optimized weight matrix}
		\EndIf
		\State Update the solution:
		\[
		x_{k+1} = x_k - g_k(x_k)
		\]
		\If{$||x_{k+1}-x_k|| < \varepsilon$}
		\State Terminate the algorithm.
		\EndIf
		\EndFor
		\State \textbf{Output:} Optimal solution $x^*$.
	\end{algorithmic}
\end{algorithm}
\newpage
It is important to note that this update procedure does not affect the convergence of OCA. Given that $ \lambda_k $ is non-increasing, it follows that the spectral radius $ \rho((R_k + H(x_k))^{-1}R_k) $ is also non-increasing. As $ \lambda_k $ decreases, the convergence rate of OCA improves. For a detailed proof, the reader is referred to \cite{zhangzgkx}.

\section{Result}
In this section, we apply the OCA and the L-M algorithm to the collected  datasets and compare the convergence speed and accuracy of the two methods based on the experimental results.
The algorithms and data used in this study were implemented in MATLAB R2023b, with the integration of the third-party CasADi\cite{andersson2019casadi} software tool.
\subsection{Real datasets experiment}
We selected three real-world cryo-electron microscopy datasets, Centriole\footnote{http://bio3d.colorado.edu/imod/files/tutorialData-1K.tar.gz}, VEEV\footnote{https://doi.org/10.5281/zenodo.11172321} and Vibrio \footnote{https://doi.org/10.5281/zenodo.11172858} to evaluate our algorithms. 
The Centriole dataset consists of a tilted sequence of 64 projections, with the projected images tilted between -61° and +65° at 2° intervals. Each projection is a 1024x1024-pixel image, with each pixel corresponding to 1.01 nm. The VEEV dataset consists of a tilted sequence of 21 projections, with the projected images tilted between -50° and +50° at 5° intervals. Each projection is a 1536x2048-pixel image, with each pixel corresponding to 0.2 nm. The Vibrio dataset\cite{10.1371/journal.pone.0215531} consists of a tilted sequence of 41 projections. The tomographic reconstructions of these samples were performed using the incline program from the MarkerAuto software suite\cite{autoHAN2015403}.

Table 1 presents the experimental results of the real datasets, indicating that the Optimal Control algorithm has a faster convergence speed compared to the L-M algorithm.

\begin{table}[h]
	\centering
	\caption{Real datasets}
	\begin{tblr}{
			colspec = {c c c c c c c},
			vlines,                     
			hlines,        
			row{1} = {font=\bfseries},            
			cell{2,4,6,8}{1} = {r=2}{c}, 
			cell{2,4,6,8}{5} = {r=2}{c},
		}
		Data & Method & $R/\mu_0$ & Iteration & Initial residual & Final residual & Time \\
		Centriole & OCA & $R=0.08I$ &4 & 1.730 & 1.196 & 0.068s \\
		& L-M & $\mu_0=0.1$& 9 &       & 1.196 & 0.193s  \\
		VEEV      & OCA & $R=17I$ & 28 & 10.413 & 1.202 & 0.176s  \\
		& L-M & $\mu_0=0.1$& 236 &      & 1.202 & 1.776s  \\
		Vibrio    & OCA & $R=0.01I$ &5 & 10.431 & 1.124 & 0.372s \\
		& L-M & $\mu_0=0.1$&7 &        & 1.124 &1.920s \\
		Adhesion belt & OCA & $R=5I$ &12 & 1.065 &0.477 &0.751s \\
		& L-M & $\mu_0=1$&6 &       &0.477 &1.357s \\
	\end{tblr}
\end{table}
\begin{figure}[htp]
	\includegraphics[scale=0.33]{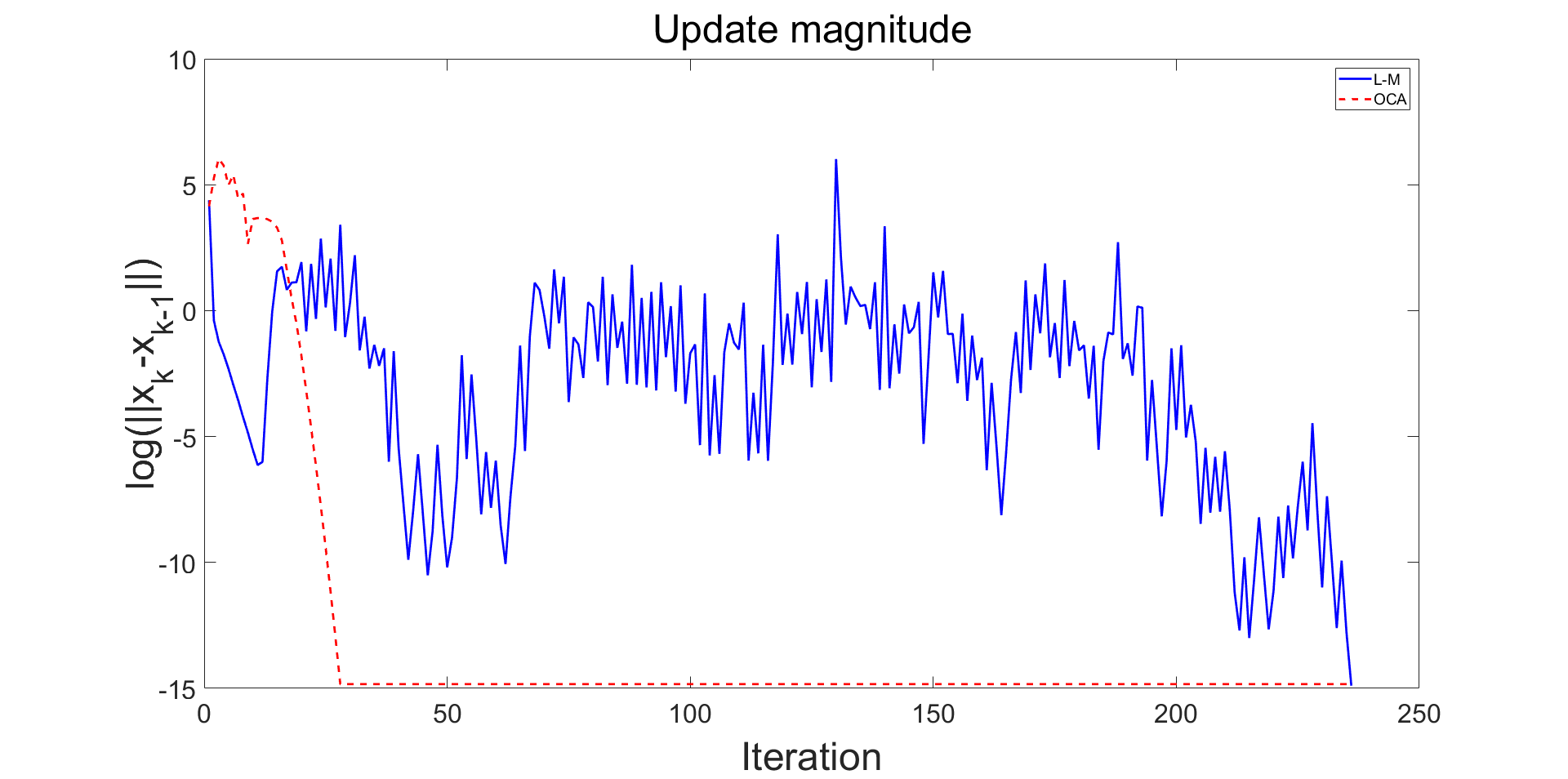}
	\caption{Convergence error}
\end{figure}
Using the VEEV dataset as an illustrative example, Figure 1 presents the error in the parameters to be optimized at each iteration, represented as $ x_k - x_{k-1} $. Notably, when $ k = 28 $, the OCA satisfied the convergence criterion of $ \varepsilon < 10^{-6} $. In contrast, the L-M algorithm exhibited oscillatory behavior, which led to an increased number of iterations. Consequently, it only met the convergence criterion at $ k = 236 $.

\subsection{Simulated datasets experiment}
To obtain a broader range of experimental results, we randomly constructed $n$ 3D points and $m$ projection images at various angles for simulated dataset experiments, where $ n = 20, 40, 60 $ and $ m = 21, 41, 64$. The projection images are sized $ 1024 \times 1024 $ pixels, and the 3D points are distributed within an $ 800 \times 800 \times 400 $ space.  To enhance the realism of the simulated dataset, we processed the constructed data as follows:
\begin{itemize}
	\item \textbf{Step 1}\quad Randomly select 5\%–30\% of the 2D projection points as invisible.
	\item  \textbf{Step 2}\quad Add Gaussian noise with a standard deviation of $ a\% \times 1024 $ to the 2D projection points, where $ a = 0.2, 2, 10 $.
	\item  \textbf{Step 3}\quad Add $ b\% $ noise to the camera parameters,  with noise magnitude based on the average of each parameter (avg), uniformly distributed in the range $[- \text{avg}, \text{avg}]$, where $ b = 5, 10 $.  
	\item  \textbf{Step 4}\quad Select 5\% of the 2D projection points as outliers and add Gaussian noise with a mean of $ 0.01 \times 1024 $ and a standard deviation of $ 0.04 \times 1024 $.
\end{itemize}

For convenience, we  denote the dataset with parameters $m=21, n=20, a=0.2, b=5$ as “m21 n20 a0.2 b5”. To ensure a fair comparison, the weight matrix $R$ in the OCA and the damping factor $\mu$ in the L-M algorithm were iteratively adjusted to promote faster convergence. The results of this experiment are presented in Table 2. 
\begin{table}[htp]
	\centering
	\caption{Simulated datasets for small noise}
	\begin{tblr}{
			width = \linewidth,    
			colspec = {c c c c c c c}, 
			vlines,                
			hlines,                
			row{1} = {font=\bfseries}, 
			cell{2,4,6,8,10,12,14,16}{1} = {c=1}{}, 
			cell{2,4,6,8,10,12,14,16}{5} = {c=1}{}, 
		}
		Data              & Method &  $ R/\mu_0 $ & Iteration & Initial residual & Final residual & Time      \\
		\SetCell[r=2]{} m21 n20 a0.2 b5 
		& OCA    & \( R=0.25I \) & 6         & \SetCell[r=2]{} 11.178 
		& 2.051          & 0.024s \\
		& L-M    & \( \mu_0=0.1 \) & 8       &                   & 2.051          & 0.032s \\
		\SetCell[r=2]{} m21 n20 a2 b5 
		& OCA    & \( R=0.25I \) & 6         & \SetCell[r=2]{} 14.683 
		& 12.378         & 0.023s \\
		& L-M    & \( \mu_0=1 \)   & 24      &                   & 12.378         & 0.090s \\
		\SetCell[r=2]{} m41 n20 a0.2 b5 
		& OCA    & \( R=0.25I \) & 5         & \SetCell[r=2]{} 8.412  
		& 1.920          & 0.034s \\
		& L-M    & \( \mu_0=0.1 \) & 6       &                   & 1.920          & 0.050s \\
		\SetCell[r=2]{} m41 n20 a2 b5 
		& OCA    & \( R=0.625I \) & 5         & \SetCell[r=2]{} 15.658 
		& 12.107         & 0.040s \\
		& L-M    & \( \mu_0=10 \)  & 17      &                   & 12.107         & 0.144s \\
		\SetCell[r=2]{} m21 n40 a0.2 b5 
		& OCA    & \( R=0.25I \) & 6         & \SetCell[r=2]{} 8.085  
		& 1.992          & 0.035s \\
		& L-M    & \( \mu_0=0.1 \) & 6       &                   & 1.992          & 0.042s \\
		\SetCell[r=2]{} m21 n40 a2 b5 
		& OCA    & \( R=I \)      & 8         & \SetCell[r=2]{} 14.192 
		& 11.834         & 0.042s \\
		& L-M    & \( \mu_0=1 \)   & 16      &                   & 11.834         & 0.010s \\
		\SetCell[r=2]{} m41 n40 a0.2 b5 
		& OCA    & \( R=0.25I \) & 5         & \SetCell[r=2]{} 7.968  
		& 1.861          & 0.058s \\
		& L-M    & \( \mu_0=0.1 \) & 6       &                   & 1.861          & 0.089s \\
		\SetCell[r=2]{} m41 n40 a2 b5 
		& OCA    & \( R=0.5I \)   & 5         & \SetCell[r=2]{} 16.647 
		& 13.512         & 0.050s \\
		& L-M    & \( \mu_0=1 \)   & 9       &                   & 13.512         & 0.153s \\
	\end{tblr}
\end{table}

When the high noise is introduced into the simulated dataset, the initial estimates of the camera parameters and the 2D point coordinates exhibit significant deviations. This suggests that the initial parameters are considerably distant from the extremum, necessitating a greater number of iterations for both the OCA and L-M algorithm. It is worth noting that, in the L-M algorithm, the damping factor \(\mu\) is continuously updated, whereas the weight matrix \(R\) in the OCA remains constant throughout the iterations. If we allow the matrix $R$ in the OCA to be updated gradually similar to \(\mu\) in the L-M algorithm, the OCA could potentially achieve faster convergence. The detailed experimental results are displayed in Table 3.

\begin{table}[htp]
	\centering
	\caption{Simulated datasets for large noise}
	\begin{tblr}{
			width = \linewidth,    
			colspec = {c c c c c c c}, 
			vlines,                
			hlines,                
			row{1} = {font=\bfseries}, 
			cell{2,4,6,8}{1} = {c=1}{}, 
			cell{2,4,6,8}{5} = {c=1}{}, 
		}
		Data              & Method & \( R_0/\mu_0 \) & Iteration & Initial residual & Final residual & Time       \\
		\SetCell[r=2]{} m21 n20 a10 b10 
		& OCA    & \( R_0=10^5I \) & 8         & \SetCell[r=2]{} 75.678 
		& 67.391         & 0.063s  \\
		& L-M    & \( \mu_0=0.01 \) & 108      &                   & 67.391         & 0.411s  \\
		\SetCell[r=2]{} m41 n20 a10 b10 
		& OCA    & \( R_0=2{\times}10^5I \) & 9 & \SetCell[r=2]{} 68.692 
		& 58.815         & 0.117s  \\
		& L-M    & \( \mu_0=0.1 \)   & 146     &                   & 58.815         & 0.100s  \\
		\SetCell[r=2]{} m21 n40 a10 b10 
		& OCA    & \( R_0=75I \)     & 8        & \SetCell[r=2]{} 70.903 
		& 65.472         & 0.080s  \\
		& L-M    & \( \mu_0=100 \)   & 168     &                   & 65.472         & 1.133s  \\
		\SetCell[r=2]{} m41 n40 a10 b10 
		& OCA    & \( R_0=10I \)     & 8        & \SetCell[r=2]{} 63.125 
		& 57.466         & 0.237s  \\
		& L-M    & \( \mu_0=0.1 \)   & 30      &                   & 57.466         & 0.452s  \\
	\end{tblr}
\end{table}
\section{Discussion}
In this paper, we introduce a novel algorithm, OCA, to optimize the BA model for the image sequence alignment of cryo-ET . We extended the OCA from a purely theoretical framework to practical applications. Through experiments on real-world datasets, the OCA exhibited a notably faster convergence rate compared to the widely adopted L-M algorithm, underscoring its practical significance. Furthermore, in experiments with simulated datasets, the advantages of the OCA became even more pronounced as the noise level increased. This indicates that, in real-world scenarios with poor initial estimates, the OCA may offer substantial benefits.

The BA model used for cryo-ET image sequence alignment involves parameters of relatively low magnitude, with matrices of only a few hundred dimensions. However, the successful application of the OCA in this context suggests its broader potential within the field of computer vision. While the implicit iterative formula of the OCA theoretically guarantees the optimal solution, as described in (\ref{impite}), it is technically challenging to implement. At present, several explicit iterative formulas have been proposed \cite{zhangzgkx}, and many explicit solutions remain to be explored. This presents a promising avenue for future research.

	\bibliography{makerbased}
\end{document}